\documentclass[11pt]{amsart}
\usepackage{graphicx,amscd,color,amsmath,amsfonts,amssymb,geometry}
\newtheorem{theorem}{Theorem}[section]

\numberwithin{equation}{section}
\begin{document}
\title[On the optimal multilinear Bohnenblust--Hille constants]{On the optimal multilinear Bohnenblust--Hille constants}
\author[D. Nu\~{n}ez et al.]{D. Nu\~{n}ez-Alarc\'{o}n\textsuperscript{*}, D. Pellegrino\textsuperscript{**} \and J.B. Seoane-Sep\'{u}lveda\textsuperscript{***} \and D. M. Serrano-Rodr\'{\i}guez\textsuperscript{*}}
\address{Departamento de Matem\'{a}tica,\newline\indent Universidade Federal da Para\'{\i}ba,\newline\indent 58.051-900 - Jo\~{a}o Pessoa, Brazil.}
\email{danielnunezal@gmail.com}
\thanks{\textsuperscript{*}Supported by Capes.}
\address{Departamento de Matem\'{a}tica,\newline\indent Universidade Federal da Para\'{\i}ba,\newline\indent 58.051-900 - Jo\~{a}o Pessoa, Brazil.}
\email{dmpellegrino@gmail.com and pellegrino@pq.cnpq.br}
\thanks{\textsuperscript{**}Supported by CNPq Grant 301237/2009-3.}
\address{Departamento de An\'{a}lisis Matem\'{a}tico,\newline\indent Facultad de Ciencias Matem\'{a}ticas, \newline\indent Plaza de Ciencias 3, \newline\indent Universidad Complutense de Madrid,\newline\indent Madrid, 28040, Spain.}
\email{jseoane@mat.ucm.es}
\thanks{\textsuperscript{***}Supported by grant MTM2012-34341.}
\address{Departamento de Matem\'{a}tica,\newline\indent Universidade Federal da Para\'{\i}ba,\newline\indent 58.051-900 - Jo\~{a}o Pessoa, Brazil.}
\email{dmserrano0@gmail.com}

\begin{abstract}
The upper estimates for the optimal constants of the multilinear
Bohnenblust--Hille inequality obtained in \cite{jfa} are here improved to:
\vspace{0.1cm}

\begin{enumerate}
\item For real scalars: $K_{n}\leq\sqrt{2}\left(  n-1\right)  ^{0.526322}$.
\item For complex scalars: $K_{n}\leq\frac{2}{\sqrt{\pi}}\left(  n-1\right)
^{0.304975}$.
\end{enumerate}
\vspace{0.1cm}
\noindent We also obtain sharper estimates for higher values of $n$. For instance,
\[
K_{n}<1.30379\left(  n-1\right)  ^{0.526322}
\]
for real scalars and $n>2^{8}$ and
\[
K_{n}<0.99137\left(  n-1\right)  ^{0.304975}
\]
for complex scalars and $n >2^{15}.$

\end{abstract}
\maketitle

\section{Preliminaries. First estimates}

Let $\mathbb{K}$ be the real or complex scalar field. The multilinear
Bohnenblust--Hille inequality asserts that for every positive integer $n\geq1$
there exists a sequence of positive scalars $\left(  C_{n}\right)
_{n=1}^{\infty}$ in $[1,\infty)$ such that
\[
\left(  \sum\limits_{i_{1},\ldots,i_{n}=1}^{N}\left\vert U(e_{i_{^{1}}}%
,\ldots,e_{i_{n}})\right\vert ^{\frac{2n}{n+1}}\right)  ^{\frac{n+1}{2n}}\leq
C_{n}\sup_{z_{1},\ldots,z_{n}\in\mathbb{D}^{N}}\left\vert U(z_{1},\ldots
,z_{n})\right\vert
\]
for all $n$-linear forms $U:\mathbb{K}^{N}\times\cdots\times\mathbb{K}%
^{N}\rightarrow\mathbb{K}$ and every positive integer $N$, where $\left(
e_{i}\right)  _{i=1}^{N}$ denotes the canonical basis of $\mathbb{K}^{N}$ and
$\mathbb{D}^{N}$ represents the open unit polydisc in $\mathbb{K}^{N}$. In the
last years a considerable effort was spent in the searching of optimal values
for the constants $C_{n}$; for details and the state-of-the-art we refer to
\cite{jfa}. The notation and terminology used in this note are the same as
those from \cite{jfa}, where it is proved that the optimal multilinear
Bohnenblust--Hille constants $\left(  K_{n}\right)  _{n=2}^{\infty}$ satisfy
\begin{equation}
K_{n}<1.65\left(  n-1\right)  ^{0.526322}+0.13\text{ (real scalars)}
\label{rea}%
\end{equation}
and
\[
K_{n}<1.41\left(  n-1\right)  ^{0.304975}-0.04\text{ (complex scalars).}%
\]
The proof of the above estimates is achieved by following a series of
technical steps. In the case of real scalars, using some previous lemmata, it
is observed that the sequence
\[
M_{n}=\left\{
\begin{array}
[c]{ll}%
\left(  \sqrt{2}\right)  ^{n-1} & \text{ if }n=1,2\\
DM_{\frac{n}{2}} & \text{ if }n\text{ is even, and}\\
DM_{\frac{n+1}{2}} & \text{ if }n\text{ is odd}%
\end{array}
\right.
\]
satisfies the multilinear Bohnenblust--Hille inequality, where $D=\frac
{e^{1-\frac{1}{2}\gamma}}{\sqrt{2}}.$ Then, using a \textquotedblleft uniform
approximation\textquotedblright\ argument, the estimate (\ref{rea}) is
achieved. In this section we remark that this final step of the proof, i.e.,
the uniform approximation argument, can be dropped and a quite simple argument
provides even better constants. In fact, from \cite{jfa} we know that, for all
$k\geq1$ and $n\geq2$, we have
\[
M_{n}=\sqrt{2}D^{k-1}\text{ whenever }n\in B_{k}=\{2^{k-1}+1,\ldots,2^{k}\}.
\]
Thus, $k-1\leq\log_{2}\left(  n-1\right)  $ and, hence,
\[
M_{n}\leq\sqrt{2}D^{\log_{2}\left(  n-1\right)  }=\sqrt{2}\left(  n-1\right)
^{\log_{2}\left(  \frac{e^{1-\frac{1}{2}\gamma}}{\sqrt{2}}\right)  }\leq
\sqrt{2}\left(  n-1\right)  ^{0.526322}.
\]
Using a similar argument (for complex scalars) it follows that
\[
M_{n}\leq\frac{2}{\sqrt{\pi}}\left(  n-1\right)  ^{\log_{2}\left(  e^{\frac{1}{2}-\frac{1}{2}\gamma}\right)  }\leq\frac{2}{\sqrt{\pi}}\left(  n-1\right)^{0.304975}
\]
for the complex scalar field. Of course, the other estimates of \cite{jfa} related to the above results can be straightforwardly improved by using these
new estimates.

\section{Sharper estimates for big values of $n$}

In this section we improve our previous estimates for {\em large} values of $n.$

\subsection{Real case}

If $\left(  C_{n}\right)  _{n=1}^{\infty}$ denotes the sequence in \cite[(4.3)]{jfa}, if we fix any $k_{0},$ it is obvious that

\[
J_{n}=\left\{
\begin{array}{ll}
C_{n} & \text{ if }n\leq2^{k_{0}},\\
DJ_{\frac{n}{2}} & \text{ if }n>2^{k_{0}}\text{ is even, and }\\
D\left(J_{\frac{n-1}{2}}\right)  ^{\frac{n-1}{2n}}\left(  J_{\frac{n+1}{2}}\right)^{\frac{n+1}{2n}} &\text{ if }n>2^{k_{0}}\text{ is odd}
\end{array}
\right.
\]
with $D=\frac{e^{1-\frac{1}{2}\gamma}}{\sqrt{2}}$, satisfies the multilinear Bohnenblust--Hille inequality. For $n>2^{k_{0}}$, let $k_{1}>k_{0}$ be such that
\[
2^{k_{1}-1}+1\leq n\leq2^{k_{1}}.
\]
Then%
\[
k_{1}-k_{0}\leq\log_{2}\left(  \frac{n-1}{2^{k_{0}-1}}\right)
\]
and
\begin{align*}
K_{n}  &  \leq J_{2^{k_{1}}}  =D^{k_{1}-k_{0}}C_{2^{k_{0}}}\\
&  \leq C_{2^{k_{0}}}D^{\log_{2}\left(  \frac{n-1}{2^{k_{0}-1}}\right)  }\\
&  =\frac{C_{2^{k_{0}}}}{D^{k_{0}-1}}\left(  n-1\right)  ^{\log_{2}D}.
\end{align*}
We thus have
\[
K_{n}\leq\frac{C_{2^{k_{0}}}}{D^{k_{0}-1}}\left(  n-1\right)  ^{\log
_{2}\left(  \frac{e^{1-\frac{1}{2}\gamma}}{\sqrt{2}}\right)  }.
\]
From \cite[Theorem 3.1]{diana} we know that
\begin{equation}
C_{2^{k_{0}}}\leq4D^{k_{0}-4} \label{gt}
\end{equation}
whenever $k_{0}\geq4.$ Thus,
\[
K_{n}\leq\frac{4}{\left(  \frac{e^{1-\frac{1}{2}\gamma}}{\sqrt{2}}\right)^{3}}\left(  n-1\right)  ^{\log_{2}\left(  \frac{e^{1-\frac{1}{2}\gamma}
}{\sqrt{2}}\right)  }<1.338887\left(  n-1\right)  ^{0.526322}.
\]
Summarizing, we have:

\begin{theorem}
If $n>16$, then
\[
K_{n}\leq\frac{4}{\left(  \frac{e^{1-\frac{1}{2}\gamma}}{\sqrt{2}}\right)^{3}}\left(  n-1\right)  ^{\log_{2}\left(  \frac{e^{1-\frac{1}{2}\gamma}
}{\sqrt{2}}\right)  }.
\]
Numerically,
\begin{equation}
K_{n}<1.338887\left(  n-1\right)  ^{0.526322}. \label{tp}%
\end{equation}

\end{theorem}

If we use the exact value of $C_{2^{k_{0}}}$ instead of estimate (\ref{gt}) we can improve (\ref{tp}) as $n$ grows. For example,
\begin{align*}
n  &  >2^{6}\Rightarrow K_{n}<1.310883\left(  n-1\right)  ^{0.526322}\\
n  &  >2^{7}\Rightarrow K_{n}<1.306156\left(  n-1\right)  ^{0.526322}\\
n  &  >2^{8}\Rightarrow K_{n}<1.303787\left(  n-1\right)  ^{0.526322}.
\end{align*}

\subsection{Complex case}

Let $\left(C_{n}\right)_{n=1}^{\infty}$ denote the sequence in \cite[Theorem 2.3]{jfa2}. If we fix any $k_{0}$, and as the authors did in \cite{jfa}, we can show that
\[
J_{n}=\left\{\begin{array}{ll}
C_{n} & \text{ if }n\leq2^{k_{0}}, \\
DJ_{\frac{n}{2}} & \text{ if }n>2^{k_{0}}\text{ is even, and}\\
D\left(  J_{\frac{n-1}{2}}\right)  ^{\frac{n-1}{2n}}\left(  J_{\frac{n+1}{2}}\right)  ^{\frac{n+1}{2n}} & \text{ if }n>2^{k_{0}}\text{ is odd,}
\end{array}
\right.
\]
with $D=e^{\frac{1}{2}-\frac{1}{2}\gamma}$, satisfies the multilinear Bohnenblust--Hille inequality. For $n>2^{k_{0}}$, by mimicking the real
case we obtain
\[
K_{n}\leq\frac{C_{2^{k_{0}}}}{D^{k_{0}-1}}\left(  n-1\right)  ^{\log_{2}\left(  e^{\frac{1}{2}-\frac{1}{2}\gamma}\right)  }.
\]
Thus, using the values of $C_{2^{k}}$ from \cite{jfa2} we have
\begin{align*}
n &  >2^{3}\Rightarrow K_{n}<1.02960973695\left(  n-1\right)  ^{0.304975}, \\
n &  >2^{4}\Rightarrow K_{n}<1.01089344604\left(  n-1\right)^{0.304975}, \\
n &  >2^{5}\Rightarrow K_{n}<1.00123230777\left(  n-1\right)  ^{0.304975}, \\
n &  >2^{6}\Rightarrow K_{n}<0.99632125476\left(  n-1\right)  ^{0.304975}, \, \ldots \\
n &  >2^{14}\Rightarrow K_{n}<0.99137409768\left(  n-1\right)  ^{0.304975}, \\
n &  >2^{15}\Rightarrow K_{n}<0.99136434217\left(  n-1\right)  ^{0.304975}, \, \ldots \\
n &  >2^{25}\Rightarrow K_{n}<0.99135459597\left(  n-1\right)  ^{0.304975}, \, \ldots \\
n &  >2^{50}\Rightarrow K_{n}<0.99135458644\left(  n-1\right)  ^{0.304975}.\\
\end{align*}

\end{document}